\documentclass[12pt, a4paper, twoside]{article}
\usepackage{amssymb}
\usepackage[leqno]{amsmath}
\usepackage{latexsym}
\usepackage{enumerate}
\usepackage{amsmath, amssymb, amscd}
\Roman{equation}

\DeclareMathOperator{\gr}{gr}
\DeclareMathOperator{\rk}{rk}
\DeclareMathOperator{\Pic}{Pic}
\DeclareMathOperator{\lcm}{lcm}

\begin{document}

\title{\textbf{\Large{Unstable Lazarsfeld-Mukai bundles of rank 2 on a certain K3 surface of Picard number 2}}}

\author{Kenta Watanabe \thanks{Nihon University, College of Science and Technology,   7-24-1 Narashinodai Funabashi city Chiba 274-8501 Japan , {\it E-mail address:watanabe.kenta@nihon-u.ac.jp}, Telephone numbers: 090-9777-1974} }

\date{}

\maketitle 


\begin{abstract}

\noindent Let $g$ and $c$ be any integers satisfying $g\geq3$ and $0\leq c\leq \lfloor\frac{g-1}{2}\rfloor$. It is known that there exists a polarized K3 surface $(X,H)$ such that $X$ is a K3 surface of Picard number 2, and $H$ is a very ample line bundle on $X$ of sectional genus $g$ and Clifford index $c$, by Johnsen and Knutsen([J-K] and [Kn]). In this paper, we give a necessary and sufficient condition for a Lazarsfeld-Mukai bundles of rank 2 associated with a smooth curve $C$ belonging to the linear system $|H|$ and a base point free pencil on $C$ not to be $H$-slope stable.

\end{abstract}

\section{Introduction} Previously, the stability of vector bundles with respect to a given polarization and the moduli spaces of stable bundles on an algebraic variety have been studied by many people. In particular, the geometry of the moduli spaces of stable vector bundles on K3 surfaces is closely connected with the Brill-Noether theory of curves via Lazarsfeld-Mukai bundles on K3 surfaces, and plays an important role in the algebraic curve theory. Recently, stable Lazarsfeld-Mukai bundles and the restrictions to curves of them have been often used to construct counterexamples to the Mercat conjecture in the Brill-Noether theory of higher rank (for example, see [S], [FO] and so on). 

Let $X$ be a K3 surface, $C$ be a smooth non-hyperelliptic curve of genus $g$ on it, and $Z$ be a divisor of degree $d$ on $C$ which forms a base point free linear system $g_d^r$. Then the Lazarsfeld-Mukai bundle $E_{C,Z}$ of rank $r+1$ associated with them is defined as the dual of the kernel $F_{C,Z}$ of the evaluation map
$${\rm{ev}}:H^0(\mathcal{O}_C(Z))\otimes\mathcal{O}_X\longrightarrow \mathcal{O}_C(Z).$$
It is well known that if $\rho(g,r,d)=g-(r+1)(g-d+r)<0$, then $E_{C,Z}$ is not simple. In particular, if $r=1$, that is, the rank of $E_{C,Z}$ is two, then the following result is known by Donagi and Morrison.

\newtheorem{thm}{Theorem}[section]

\begin{thm} {\rm{([D-M], Lemma 4.4)}}. Let the notations be as above. If $E_{C,Z}$ is of rank 2 and not simple, then there exist two line bundles $M$ and $N$ on $X$ and a 0-dimensional subscheme $Z^{'}\subset X$ of finite length such that

$\;$

{\rm{(a)}} $h^0(M)\geq2,\;h^0(N)\geq2;$ 

\smallskip

\smallskip

{\rm{(b)}} $N$ is base point free;

\smallskip

\smallskip

{\rm{(c)}} There exists an exact sequence
$$0\longrightarrow M\longrightarrow E_{C,Z}\longrightarrow N\otimes\mathcal{J}_{Z^{'}}\longrightarrow0.$$
Moreover, if $h^0(M\otimes N^{\vee})=0$, then the length of $Z^{'}$ is zero.

\end{thm}

\noindent The extension of the rank one torsion-free sheaf as in Theorem 1.1 is called the Donagi-Morrison's extension. Moreover, if $h^0(M\otimes N^{\vee})=0$, then $E_{C,Z}\cong N\oplus M$. This means that if $\rho(g,1,d)=2d-g-2<0$ and $C$ is very ample, the Lazarsfeld-Mukai bundle $E_{C,Z}$ of rank 2 associated with $C$ and the base point free pencil $|Z|=g_d^1$ on $C$ is not $C$-slope stable. Therefore, conversely, it is natural and interesting to consider the problem of when $E_{C,Z}$ is not $C$-slope stable, in the case where $\rho(g,1,d)\geq0$. Recently, the following work concerning with the description of the moduli space of unstable Lazarsfeld-Mukai bundles of rank 2 is known.

\begin{thm} {\rm{([LC], Theorem 4.3)}}. Let $X$ be a K3 surface, and $H$ be a very ample line bundle on $X$. Assume that the general smooth curves in $|H|$ have Clifford dimension one and maximal gonality $k=\lfloor\frac{g+3}{2}\rfloor$, then:

\smallskip

\smallskip

\noindent {\rm{(i)}} If $\rho(g,1,d)>0$, then any dominating component of $\mathcal{W}_d^1(|H|)$ corresponds to $H$-slope stable Lazarsfeld-Mukai bundles associated with smooth curves in $|H|$ and base point free pencils on them. In particular, if $C\in |H|$ is a general smooth curve, then the variety $W_d^1(C)$ which is the fiber over $C$ of the natural projection $\pi: \mathcal{W}_d^1(|H|)\rightarrow |H|$ is reduced and has the expected dimension $\rho(g,1,d)$.

\smallskip

\smallskip

\noindent {\rm{(ii)}} If $\rho(g,1,k)=0$ and $C\in |H|$ is a general smooth curve, then $W_k^1(C)$ has dimension 0.

\end{thm}

\noindent The assertion of (i) in the above theorem means that if $C\in |H|$ is general and $\rho(g,1,d)>0$, then the Lazarsfeld-Mukai bundle associated with a general complete base point free pencil $g_d^1$ on $C$ is $H$-slope stable. Moreover, the assumption that $H$ has maximal Clifford index $\lfloor\frac{g-1}{2}\rfloor$ is essentially used to prove it. However, this theorem does not answer the above question. In this paper, we will give a necessary and sufficient condition for $E_{C,Z}$ on a certain K3 surface whose Picard lattice contains $H=\mathcal{O}_X(C)$ as a generator of it not to be $C$-slope stable. Our main theorem is as follows. 

\begin{thm}. Let $d$ and $g$ be integers with $3\leq d\leq \lfloor\frac{g+3}{2}\rfloor$. Let $X$ be a K3 surface with $Pic(X)=\mathbb{Z}H\oplus\mathbb{Z}F$ such that $H$ is a base point free line bundle on $X$ with $H^2=2g-2$, and $F$ is an elliptic pencil on $X$ satisfying $H.F=d$. Then $H$ is very ample, and if $C\in |H|$ is a smooth curve and $Z$ is a divisor on $C$ which forms a base point free pencil on $C$, then the Lazarsfeld-Mukai bundle $E_{C,Z}$ associated with them has the following properties.

\smallskip

\smallskip

\noindent {\rm{(i)}} If $|K_C\otimes\mathcal{O}_C(-Z)|=\emptyset$, then $E_{C,Z}$ is $H$-slope stable.

\smallskip

\smallskip

\noindent {\rm{(ii)}} If $E_{C,Z}$ is not $H$-slope stable, then $\mathcal{O}_C(Z)=F|_C$ or $H\otimes F^{\vee}|_C$, and the latter case occurs precisely when $d=g-1$.

\end{thm}

\noindent We can easily see that the consequence of the above theorem that $\mathcal{O}_C(Z)=F|_C$ (resp. $\mathcal{O}_C(Z)=H\otimes F^{\vee}|_C$) means that $H\otimes F^{\vee}\subset E_{C,Z}$ (resp. $F\subset E_{C,Z}$). Hence, if $d\leq g-1$, it gives a necessary and sufficient condition for $E_{C,Z}$ not to be $H$-slope stable.

In our previous works, it is known that $E_{C,Z}$ as in Theorem 1.3 is ACM and initialized with respect to $H$ ([W1], Proposition 2.3). Moreover, we can easily see that $F$ and $H\otimes F^{\vee}$ are unique ACM and initialized line bundles with respect to $H$ in the above theorem. Hence, in order to prove Theorem 1.3, we show that any saturated sub-line bundle $L$ of $E_{C,Z}$ with $H.L\geq g-1$ is ACM with respect to $H$.

Our plan of this paper is as follows. In section 2, we recall the notion of Mumford-Takemoto stability for vector bundles on projective varieties, and give a necessary condition for Lazarsfeld-Mukai bundles of rank 2 on K3 surfaces not to be stable. In section 3, we recall the existence theorem of the polarized K3 surface $(X,H)$ as in Theorem 1.3, and prove our main theorem.

$\;$

\noindent {\bf{Notations and conventions}}. We work over the complex number field $\mathbb{C}$. A curve and a surface are smooth projective. For a curve or a surface $Y$, we denote by $K_Y$ the canonical line bundle of $Y$ and denote by $|D|$ the linear system defined by a divisor $D$ on $Y$. 

For a surface $X$, $\Pic(X)$ denotes the Picard lattice of $X$. A regular surface $X$ (i.e., a surface $X$ with $h^1(\mathcal{O}_X)=0$) is called a K3 surface if the canonical bundle $K_X$ of it is trivial. For a vector bundle $E$, we denote by $E^{\vee}$ the dual of it, and denote by $\rk E$ the rank of $E$. If we fix a very ample line bundle $H$ as a polarization on $X$, then we write $E\otimes H^{\otimes l}=E(l)$ for a vector bundle $E$ on $X$. We will say that a vector bundle $E$ is initialized with respect to a given polarization $H$, if 
$$H^0(E)\neq0,\text{ and }H^0(E(-1))=0.$$
%We will say that a vector bundle $E$ is ACM with respect to a given polarization $H$ on a regular surface $X$, if it satisfies $H^1(E(l))=0$, for any integer $l\in\mathbb{Z}$.

\section{Slope stability of Lazarsfeld-Mukai bundles on K3 surfaces}

In this section, we recall the definition of Mumford-Takemoto stability (i.e., slope stability) of vector bundles with respect to a given polarization on a projective variety, and prepare some propositions about the slope stability of Lazarsfeld-Mukai bundles of rank 2 on K3 surfaces to prove our main theorem.

\newtheorem{df}{Definition}[section]

\begin{df} Let $X$ be a smooth projective variety, $H$ be a very ample line bundle on $X$, and let $E$ be a torsion free sheaf on $X$ of rank $r$. Then the $H$-slope of $E$ is defined as follows;
$$\mu_H(E)=\frac{c_1(E).H}{r}.$$
\noindent $E$ is called {\rm{$\mu_H$-semistable (resp. $\mu_H$-stable)}} if for any subsheaf $0\neq F\subset E$ with $\rk F<\rk E$, we have $\mu_H(F)\leq\mu_H(E)$ {\rm{(}}resp. $\mu_H(F)<\mu_H(E)${\rm{)}}.\end{df}

\noindent In this paper, we often call a $\mu_H$-stable (resp. $\mu_H$-semistable) bundle a $H$-slope stable (resp. $H$-slope semistable) bundle. It is well known that for a vector bundle $E$ on $X$, there is a unique filtration called the {\it{Harder-Narasimhan {\rm{(}}HN for short{\rm{)}} filtration}} 
$$0=E_0\subset E_1\subset\cdots\subset E_n=E,$$
\noindent such that $E_i$ is locally free and $E_i/E_{i-1}$ is a torsion free and $\mu_H$-semistable sheaf, for $1\leq i\leq n$, and $\mu_H(E_{i+1}/E_i)<\mu_H(E_i/E_{i-1})$, for $1\leq i\leq n-1$. Moreover, such a filtration satisfies the following inequality
$$\mu_H(E_1)>\mu_H(E_2)>\cdots>\mu_H(E).$$
\noindent Obviously, if $E$ is not semistable, $n\geq2$. Then the sheaf $E_1$ is called the maximal destabilizing sheaf of $E$. Moreover, if a vector bundle $E$ is $\mu_H$-semistable, there exists a filtration called a {\it{Jordan-H\"{o}lder {\rm{(}}JH for short{\rm{)}} filtration}}
$$0=JH_0(E)\subset JH_1(E)\subset\cdots\subset JH_m(E)=E,$$
\noindent such that $\gr_i(E):=JH_i(E)/JH_{i-1}(E)$ is a torsion free and $\mu_H$-stable sheaf whose slope is equal to $\mu_H(E)$ for $1\leq i\leq m$. 

From now on, we assume that $X$ is a K3 surface. First of all, we remark the following assertion.

\newtheorem{lem}{Lemma}[section]

\begin{lem}{{\rm{([LC], Lemma 3.2)}}}. Let $E$ and $Q$ be torsion free sheafs such that $\rk E\geq2$. If $E$ is globally generated off a finite number of points, $h^2(E)=0$ and there exists a surjective morphism $\varphi:E\longrightarrow Q$, then $h^0(Q^{\vee\vee})\geq2$. In particular, if the rank of $Q$ is one, then $Q^{\vee\vee}$ is a non-trivial and base point free line bundle.\end{lem}

\noindent It is well known that the Lazarsfeld-Mukai bundle $E_{C,Z}$ associated with a smooth curve $C$ on $X$ of genus $g\geq3$ and a base point free divisor $Z$ on $C$ is globally generated off the base points of $|K_C\otimes\mathcal{O}_C(-Z)|$ (cf. [A-F], Proposition 2.1). Hence, by Lemma 2.1, we find that if $\rk E_{C,Z}=2$, that is, $|Z|$ is a pencil and $|K_C\otimes\mathcal{O}_C(-Z)|\neq\emptyset$, then $(E_{C,Z}/L)^{\vee\vee}$ is a non-trivial and base point free line bundle, for any saturated sub-line bundle $L\subset E_{C,Z}$.
By using this fact, we prove the following proposition.

\newtheorem{prop}{Proposition}[section]

\begin{prop}. Let $X$ be a K3 surface, $C$ be a very ample smooth curve of genus $g$ on $X$, and $Z$ be a divisor on $C$ such that $|Z|$ is a base point free pencil. Moreover, we set $H=\mathcal{O}_X(C)$ and assume that $|K_C\otimes\mathcal{O}_C(-Z)|\neq\emptyset$. Then if a line bundle $L$ is a saturated sub-line bundle of $E_{C,Z}$ with $L.H\geq g-1$, then:

\smallskip

\smallskip

\noindent{\rm{(i)}} $L^2\geq0$. Moreover, if $L^2=0$, then $L.H=g-1$.

\smallskip

\smallskip

\noindent{\rm{(ii)}} If $(H\otimes L^{\vee})^2=0$, then:

\smallskip

\smallskip

{\rm{(a)}} $|H\otimes L^{\vee}|$ is an elliptic pencil with $H\otimes L^{\vee}|_C\cong\mathcal{O}_C(Z)$.

\smallskip

{\rm{(b)}} $E_{C,Z}/L\cong H\otimes L^{\vee}$.

\smallskip

{\rm{(c)}} $h^1(L)=0$.

\end{prop}

\noindent Before the proof of this proposition, we recall the following classical results about the classification of linear systems on K3 surfaces.

\begin{lem} {\rm{([SD], Proposition 2.7)}}. Let $L$ be a numerical effective line bundle on a K3 surface $X$. Then $|L|$ is not base point free if and only if there exists an elliptic curve $F$, a smooth rational curve $\Gamma$, and an integer $k\geq2$ such that $F.\Gamma=1$ and $L\cong\mathcal{O}_X(kF+\Gamma)$.
\end{lem}

\begin{lem} {\rm{([SD], Proposition 2.6)}}. Let $L$ be a line bundle on a K3 surface $X$ such that $|L|\neq\emptyset$. Assume that $|L|$ is base point free. Then one of the following cases occurs.

$\;$

{\rm{(i)}} $L^2>0$ and the general member of $|L|$ is a smooth irreducible curve of genus $\displaystyle\frac{L^2}{2}+1$.

\smallskip

\smallskip

{\rm{(ii)}} $L^2=0$ and $L\cong\mathcal{O}_X(kF)$, where $k\geq1$ is an integer and $F$ is a smooth curve of genus one. In this case, $h^1(L)=k-1$.

\end{lem}

\noindent We note that, a linear system $|C|$ given by an irreducible curve $C$ with $C^2>0$ is base point free. Hence, by Lemma 2.3, any line bundle $L$ on a K3 surface has no base point outside its fixed components. 

$\;$

\noindent {\bf{Proof of Proposition 2.1}}. (i) Since $L$ is saturated, $(E_{C,Z}/L)^{\vee\vee}\cong H\otimes L^{\vee}$, and $L.(H\otimes L^{\vee})\leq\deg Z$. Since $|Z|$ is a pencil such that $|K_C\otimes\mathcal{O}_C(-Z)|\neq\emptyset$, by the Riemann-Roch theorem, we have $g+1-\deg Z=h^1(\mathcal{O}_C(Z))>0.$
Therefore, we have $g-1\leq L.H\leq\deg Z+L^2\leq g+L^2.$
Hence, we have $L^2\geq0.$ Assume that $L^2=0$. Then we have $L.H=g-1$ or $g$. If $L.H=g$, we have $(H\otimes L^{\vee})^2=-2$. However, this contradicts the fact that $|H\otimes L^{\vee}|$ is base point free. Therefore, we have the assertion. 

\smallskip

\noindent (ii) By the assertion of (i), we have $h^0(L)\geq2$. Since $h^0(L^{\vee}\otimes E_{C,Z})>0$, by the exact sequence
$$0\longrightarrow H^0(\mathcal{O}_C(Z))^{\vee}\otimes L^{\vee}\longrightarrow L^{\vee}\otimes E_{C,Z}\longrightarrow L^{\vee}\otimes K_C\otimes\mathcal{O}_C(-Z)\longrightarrow 0,$$
we have $h^0(L^{\vee}\otimes K_C\otimes\mathcal{O}_C(-Z))>0$ and hence, we have 
$$\deg(L^{\vee}\otimes K_C\otimes\mathcal{O}_C(-Z))\geq0.$$
Moreover, since $(H\otimes L^{\vee})^2=0$, we have
$$\deg(L^{\vee}\otimes K_C\otimes\mathcal{O}_C(-Z))=L.(H\otimes L^{\vee})-\deg Z\leq0.$$
We have $\deg(L^{\vee}\otimes K_C\otimes\mathcal{O}_C(-Z))=0$. Hence, we have $H\otimes L^{\vee}|_C\cong\mathcal{O}_C(Z)$ and $E_{C,Z}/L\cong H\otimes L^{\vee}$. Since $|K_C\otimes\mathcal{O}_C(-Z)|\neq\emptyset$, by Lemma 2.1, $H\otimes L^{\vee}$ is base point free. Therefore, by Lemma 2.3 (ii), there exists an elliptic pencil $F$ on $X$ and an integer $r\geq 1$ such that $H\otimes L^{\vee}\cong F^{\otimes r}$. By the exact sequence
$$0\longrightarrow L^{\vee}\longrightarrow H\otimes L^{\vee}\longrightarrow \mathcal{O}_C(Z)\longrightarrow 0,$$
we have
$$r+1=h^0(H\otimes L^{\vee})\leq h^0(\mathcal{O}_C(Z))=2.$$
Hence, we have $r=1$. Therefore, $|H\otimes L^{\vee}|$ is an elliptic pencil. By using the above exact sequence again, we also have $h^1(L)=0$. $\hfill\square$

\section{Proof of Theorem 1.3}

In this section, we prove our main theorem. First of all, we recall the existence theorem for K3 surfaces as in Theorem 1.3.

\begin{thm} {\rm{([Kn], Proposition 2.1 and Lemma 2.2)}}. Let $d$ and $g$ be integers satisfying $g\geq3$ and $2\leq d\leq\lfloor\frac{g+3}{2}\rfloor$. Then there exists a K3 surface $X$ with $\Pic(X)=\mathbb{Z}H\oplus\mathbb{Z}F$ such that $H$ is a base point free line bundle with $H^2=2(g-1)$ and $|F|$ is an elliptic pencil satisfying $H.F=d$. Moreover, the Clifford index of $H$ is $d-2$. \end{thm}

\noindent In Theorem 3.1, any smooth curve $C\in |H|$ has gonality $d$, and it is computed by the restriction to $C$ of the elliptic pencil $|F|$. Hence, the Clifford dimension of any smooth curve $C\in |H|$ is one. In particular, if $d=2$ then, $|H|$ is a hyperelliptic linear system. Before the proof of Theorem 1.3, we remark the following lemmas.

\begin{lem} {\rm{([W2], Proposition 5.1)}}. Let the notations be as in Theorem 3.1, and assume that $3\leq d\leq \lfloor\frac{g+3}{2}\rfloor$. Then the following statements hold.

\smallskip

\smallskip

\noindent {\rm{(i)}} $X$ contains a $(-2)$-curve if and only if $d|g$. Moreover, in this case, if we set $g=md$, then $H\otimes F^{\vee\otimes m}$ is a unique $(-2)$-vector up to sign in $\Pic(X)$ and the member of $|H\otimes F^{\vee\otimes m}|$ is a $(-2)$-curve.

\smallskip

\smallskip

\noindent {\rm{(ii)}} Let $L$ be a line bundle which is given by an elliptic curve on $X$. Then;

\smallskip

{\rm{(a)}} If $d|g$, then $L=F$.

{\rm{(b)}} If $d\mid \hspace{-.67em}/g$, then $L=F$ or there exist integers $m$ and $n$ such that 
$$n>0,\;n(g-1)+md=0,\;n(g-1)=\lcm(g-1,d),$$
$\text{ and }L=H^{\otimes n}\otimes F^{\otimes m}$.

\smallskip

\smallskip

\noindent {\rm{(iii)}} $H$ is very ample.

\smallskip

\smallskip

\noindent {\rm{(iv)}} If $d\leq g-1$, then $H\otimes F^{\vee}$ is base point free.

\end{lem}

\noindent {\bf{Proof}}. See the proof of Proposition 5.1 in [W2].$\hfill\square$

$\;$

\begin{lem} Let the notations be as in Theorem 1.3. Let $L$ be a line bundle on $X$. If $L$ and $H\otimes L^{\vee}$ are base point free, and $h^1(L)=h^1(H\otimes L^{\vee})=0$, then $L=F$ or $H\otimes F^{\vee}$.\end{lem}

\noindent{\bf{Proof}}. Let $L=H^{\otimes s}\otimes F^{\otimes t}$. By the assumption, we have
$$h^1(F^{\otimes t})=h^1(L\otimes H^{\vee\otimes s})=h^1(H^{\otimes s}\otimes L^{\vee})=0.$$
Hence, we have $|t|\leq1$. Since $L$ and $H\otimes L^{\vee}$ are initialized with respect to $H$, we have $t\neq0$. Moreover, if $t=1$, then $s=0$, and if $t=-1$, then $s=1$. Hence, we have the assertion.$\hfill\square$

$\;$

\noindent {\bf{Proof of Theorem 1.3}}. (i) Assume that $|K_C\otimes\mathcal{O}_C(-Z)|=\emptyset$. Since $|Z|$ is a pencil, by the Riemann-Roch theorem, we have $\deg Z=g+1$. Assume that $E_{C,Z}$ is not $H$-slope stable, and let $L$ be a saturated sub-line bundle with $L.H\geq g-1$. Then we have $(E_{C,Z}/L)^{\vee\vee}\cong H\otimes L^{\vee}$ and $L.(H\otimes L^{\vee})\leq g+1$. Since
$$g-1\leq L.H\leq g+1+L^2,$$
we have $L^2\geq -2$. By Lemma 3.1 (i), if $L^2=-2$, then $d|g$. If we let $m=\displaystyle\frac{g}{d}$, then we have $L=H\otimes F^{\vee\otimes m}$. Hence, we have the contradiction
$L.H=g-2$. Therefore, we have $L^2\geq0$, and hence, we have $h^0(L)\geq2$. Since $L\subset E_{C,Z}$, we have
$$h^0(L)\leq h^0(E_{C,Z})=g+3-\deg Z=2.$$
Therefore, we have $h^0(L)=2$. By the Riemann-Roch theorem, we have $h^1(L)=0$ and $L^2=0$. By easy computation, we have $L=F$ or $L=H^{\otimes n}\otimes F^{\otimes m}$ with $n>0$ and $n(g-1)+md=0$. 
 
Assume that $L=F$. By the exact sequence
$$0\longrightarrow H^0(\mathcal{O}_C(Z))^{\vee}\otimes F^{\vee}\longrightarrow F^{\vee}\otimes E_{C,Z}\longrightarrow F^{\vee}\otimes K_C\otimes\mathcal{O}_C(-Z)\longrightarrow0,$$
we have
$$h^0(F^{\vee}\otimes K_C\otimes \mathcal{O}_C(-Z))=h^0(F^{\vee}\otimes E_{C,Z})>0.$$
Since $\deg Z=g+1$, we have 
$$g-d-3=\deg (F^{\vee}\otimes K_C\otimes \mathcal{O}_C(-Z))\geq0,$$
and hence, we have $d\leq g-3$. However, since $H.F=d$, this contradicts the assumption that $H.L\geq g-1$.

We consider the latter case. Since
$$n(g-1)=n(2g-2)+md=L.H=L.(H\otimes L^{\vee})\leq g+1,$$
and $g\geq3$, we have $n\leq2$. If $n=2$, then we have $d=g=3$. However, since $n(g-1)+md=0$, we have the contradiction $md=-4$. Hence, we have $n=1$. By the exact sequence
$$0\longrightarrow H^0(\mathcal{O}_C(Z))^{\vee}\otimes L^{\vee}\longrightarrow L^{\vee}\otimes E_{C,Z}\longrightarrow L^{\vee}\otimes K_C\otimes\mathcal{O}_C(-Z)\longrightarrow0,$$
we have 
$$h^0(L^{\vee}\otimes K_C\otimes\mathcal{O}_C(-Z))=h^0(L^{\vee}\otimes E_{C,Z})>0.$$
However, we have the contradiction
$$\deg(L^{\vee}\otimes K_C\otimes\mathcal{O}_C(-Z))=\deg(F^{\vee\otimes m}\otimes\mathcal{O}_C(-Z))=-2<0.$$
Therefore, $E_{C,Z}$ is $H$-slope stable.

$\;$

\noindent (ii) Since $E_{C,Z}$ is not $H$-slope stable, by the assertion of (i), we have $|K_C\otimes\mathcal{O}_C(-Z)|\neq\emptyset$. Let $L$ be a saturated sub-line bundle of $E_{C,Z}$ with $L.H\geq g-1$. Since $|H\otimes L^{\vee}|$ is base point free, by Lemma 2.3, if $(H\otimes L^{\vee})^2>0$, we have $h^1(H\otimes L^{\vee})=0$. If $(H\otimes L^{\vee})^2=0$, by Proposition 2.1 (ii), $|H\otimes L^{\vee}|$ is an elliptic pencil, and hence, we have $h^1(H\otimes L^{\vee})=0$. 

 First of all, we show that $L$ is base point free. Since $|K_C\otimes\mathcal{O}_C(-Z)|\neq\emptyset$, by Proposition 2.1 (i), we have $L^2\geq0$, and hence, we have $h^0(L)\geq2$. 

\noindent If $d\mid \hspace{-.67em}/g$, then, by Lemma 3.1 (i), there is no $(-2)$-curve on $X$. Hence, $L$ is base point free. Assume that $d|g$. By Lemma 3.1 (i), if we let $m=\displaystyle\frac{g}{d}$, then the member $\Gamma$ of $|H\otimes F^{\vee\otimes m}|$ is a unique $(-2)$-curve on $X$. Moreover, by Lemma 3.1 (ii), $F$ is a unique elliptic pencil on $X$, and $F.\Gamma=d\geq3$. Hence, by Lemma 2.2, it is sufficient to show that $L.\Gamma\geq0$. 

We consider the case where $L^2=0$. Since, by Proposition 2.1 (i), $L.H=g-1$, we have $(H\otimes L^{\vee})^2=0$. Hence, by Proposition 2.1 (ii) and Lemma 3.1 (ii), we have $H\otimes L^{\vee}=F$, that is, $L=H\otimes F^{\vee}$. Then we have $d\neq g$. In fact, since
$$2g-2-d=H.(H\otimes F^{\vee})=H.L\geq g-1,$$
we have $d\leq g-1$. Since $m\geq2$, we have
$$L.\Gamma=(H\otimes F^{\vee}).\Gamma=g-2-d=(m-1)d-2>0.$$

We consider the case where $L^2>0$. We assume that $L.\Gamma<0$. Then we have
$$(H\otimes L^{\vee}).\Gamma=g-2-L.\Gamma\geq g-1.$$
If we let $\Delta$ be the fixed component of $L$, then $\Delta\sim r\Gamma$ ($r\geq1$). Since $|K_C\otimes\mathcal{O}_C(-Z)|\neq\emptyset$, we have
$$g\geq (H\otimes L^{\vee}).L\geq (H\otimes L^{\vee}).\Delta=r(H\otimes L^{\vee}).\Gamma\geq r(g-1).$$
Since $g\geq3$, we have $r=1$. Here, let $L^{'}$ be the movable part of $L$. Since $L^2>0$, if $(L^{'})^2=0$, we have $L^{'}.\Gamma>1$. This contradicts the assumption that 
$L.\Gamma<0$. Since $|H\otimes L^{\vee}|$ is base point free, if $(L^{'})^2>0$, we have $(H\otimes L^{\vee})L^{'}\geq2$. However, since
$$g\geq (H\otimes L^{\vee}).L\geq (H\otimes L^{\vee}).L^{'}+g-1,$$
we have the contradiction $(H\otimes L^{\vee}).L^{'}\leq1$. Hence, $|L|$ is base point free. If $L^2>0$, then, by Lemma 2.3, we have $h^1(L)=0$. If $L^2=0$, by Proposition 2.1 (i), we have $L.H=g-1$, and hence, $(H\otimes L^{\vee})^2=0$. Therefore, by Proposition 2.1 (ii), we have $h^1(L)=0$. Hence, by Lemma 3.2, we have $L=F$ or $H\otimes F^{\vee}$. By the exact sequence
$$0\longrightarrow H^0(\mathcal{O}_C(Z))^{\vee}\otimes L^{\vee}\longrightarrow L^{\vee}\otimes E_{C,Z}\longrightarrow L^{\vee}\otimes K_C\otimes\mathcal{O}_C(-Z)\longrightarrow 0,$$
we have $h^0(L^{\vee}\otimes K_C\otimes\mathcal{O}_C(-Z))\geq h^0(L^{\vee}\otimes E_{C,Z})>0.$ Hence, we have
$$\deg(L^{\vee}\otimes K_C\otimes\mathcal{O}_C(-Z))=2g-2-L.H-\deg Z\geq0.$$
If $L=H\otimes F^{\vee}$, we have $\deg Z\leq d$. Since $d$ is the gonality of $C$, we have $\deg Z=d$, and hence, we have $F|_C\cong\mathcal{O}_C(Z)$. Assume that $L=F$. Since $d=L.H\geq g-1$, we have
$\deg Z\leq 2g-2-d\leq g-1$. Since $d=L.(H\otimes L^{\vee})\leq \deg Z$, we have $d=g-1$. Since $\deg(H\otimes F^{\vee}\otimes\mathcal{O}_C(-Z))=0$, we have $H\otimes F^{\vee}|_C\cong\mathcal{O}_C(Z)$.$\hfill\square$

$\;$

\noindent {\bf{Acknowledgements}}

\smallskip

\smallskip

\noindent The author is partially supported by Grant-in-Aid for Scientific Research (16K05101), Japan Society for the Promotion of Science.

\end{document}